\documentclass{article}
\usepackage{amsmath,amssymb,amsthm}

\usepackage{amsrefs}

\usepackage[letterpaper,right=4cm,left=4cm]{geometry}                
\usepackage{enumitem}

\newcommand{\ZZ}{\ensuremath{\mathbb Z}}

\newcommand{\FF}{\ensuremath{\mathbb F}}
\newcommand{\cA}{\ensuremath{{\cal A}}}

\newtheorem{theorem}{Theorem} 
\newtheorem*{theorem*}{Theorem}

\theoremstyle{definition}
\newtheorem{definition}[theorem]{Definition}

\DeclareMathOperator{\Bose}{\normalfont\textsc{BoseCh}}
\DeclareMathOperator{\Singer}{\normalfont\textsc{Singer}}

\title{Constructing Thick $B_h$-sets}
\author{Kevin O'Bryant\\City University of New York\\(College of Staten Island and The Graduate Center)}
\date{\today}

\begin{document}
    \maketitle
\begin{abstract}
    A subset ${\cal A}$ of a commutative semigroup $X$ is called a $B_h$-set in $X$ if the only solutions to
    \[a_1+\dots+a_h = b_1 + \cdots +b_h, \qquad a_i,b_i \in \cA\]
    are the trivial solutions $\{a_1,\dots,a_h\} = \{b_1,\dots,b_h\}$ (as multisets). With $h=2$ and $X=\ZZ$, these sets are also known as Sidon sets, Golomb Rulers, and Babcock sets. In this work, we generalize constructions of Bose-Chowla and Singer and give the resultant bounds on the diameter of a $k$ element $B_h$-set in $\ZZ$ for small $k$. We conclude with a list of open problems.
\end{abstract}

\section{Introduction}
A subset ${\cal A}$ of a commutative semigroup $X$ is called a $B_h$-set in $X$ if the only solutions to
    \[a_1+\dots+a_h = b_1 + \cdots +b_h, \qquad a_i,b_i \in \cA\]
are the trivial solutions $\{a_1,\dots,a_h\} = \{b_1,\dots,b_h\}$ (as multisets). With $h=2$ and $X=\ZZ$, these sets are also known as Sidon sets, Golomb Rulers, and Babcock sets. For an extensive bibliography of related mathematics literature we direct the reader to~\cite{bibliography}. The purpose of this work is to give new parameterized constructions of $B_h$-sets for $h\ge3$, and to give criteria on the parameters for these sets to be affinely inequivalent.

One application of $B_h$-sets in $\ZZ$ is in electrical engineering; this literature starts in Babcock~\cite{1953.Babcock} and continues for dozens of articles in IEEE journals not covered by Math Sci-Net. Specifically, a nonlinear amplifier for channel frequencies $a_1,a_2,a_3,\dots$ produces ``ghost'' signals at frequencies of the form $a_1+a_2, a_1+a_2-a_3$, and so on. The strongest relevant ghosts are at $a_1+a_2-a_3$ (third order intermodulation) and $a_1+a_2+a_3-a_4-a_5$ (fifth order intermodulation). Thus, the set of frequencies should avoid equations of the sort $a_4=a_1+a_2-a_3$ and $a_6=a_1+a_2+a_3-a_4-a_5$. That is, to avoid third order intermodulation, the channels should form a $B_2$-set, and to avoid fifth order intermodulation, the channels should be a $B_3$-set.

The first published usage of the ``$B_h$'' terminology that we have found is in the introduction of the famous Erd\H{o}s~\&~Tur\'{a}n paper~\cite{MR6197}, where they state ``Such sequences, called $B_2$ sequences by Sidon, occur in the theory of Fourier series.'' Singer~\cite{Singer} had already constructed thick finite $B_h$-sets in 1939, and Bose gave a different thick finite construction of $B_2$-sets in~\cite{1942.Bose}, which was generalized to $B_h$-sets by Bose~\&~Chowla in~\cite{BoseChowla}. The constructions given in this work subsume those of Singer and Bose~\&~Chowla.

\begin{definition}\label{def:Bose}
    For an integer $h\geq 2$ and a prime power $q$, set $M=q^h-1$. For a generator $\theta$ of the multiplicative group $\FF_{q^h}^\times$, and  $b\in \ZZ/(q^h-1)\ZZ$ for which $\theta^b$ has algebraic degree $h$ over $\FF_q$, we define the set
    \[
    \Bose_h(q,b) := \left\{a \in \ZZ/M\ZZ : \quad \theta^a = \theta^b +v, \quad v\in \FF_q\right\} .
    \]
\end{definition}

\begin{definition}\label{def:Singer}
    For an integer $h\geq 2$, a prime power $q$, set $M=\tfrac{q^{h+1}-1}{q-1}$. For a generator $\theta$ of the multiplicative group $\FF_{q^{h+1}}^\times$, and $b\in \ZZ/M\ZZ$ for which $\theta^b$ has algebraic degree $h+1$, we define the set
    \[\Singer_h(q,b) := \left\{a  \in \ZZ/M\ZZ :
    \quad\theta^a = u \theta^b +v, \quad u,v \in \FF_q\right\}.
    \]
    We comment that it may seem that the modulus should be $q^{h+1}-1$. However, for any $a$ with $\theta^a=u\theta^b+v$, one also has $\theta^{a+(q^{h+1}-1)/(q-1)}=u_1\theta^b+v_1$, where $u_1,v_1$ are in $\FF_q$ because all of $u,v,\theta^{(q^{h+1}-1)/(q-1)}$ are in $\FF_q$.
\end{definition}

A cautious reader will object that the choice of the generator impacts the right side of these definitions, and so should be included in the notation $\Bose_h(q,b)$ and $\Singer_h(q,b)$. While the choice of $\theta$ does matter, we will eventually show that it does not matter in a meaningful way. To avoid this technicality, we set $\theta$ in the above definitions to be a root of the Conway polynomial~\cite{2004.Heath.Loehr} that generates the appropriate field. The only facts about Conway polynomials that we will use is that for each prime power $q=p^e$, the Conway polynomial $C_{p,e}(\theta)\in \FF_p[\theta]$ is uniquely defined, irreducible, and
    \[\FF_p[\theta]/C_{p,e}(\theta) \cong \FF_q,\qquad \langle \theta \rangle = \FF_q^\times.\]
There are additional properties that make Conway polynomials a computationally pleasant approach to working in finite fields, particularly concerning subfields, and Luebeck~\cite{Luebeck} has provided an extensive database. The specific presentation of the finite fields is not relevant to the theory in this work, and is only useful if one wants to compare explicit computations.

The $b=1$ cases of the following theorem are exactly the constructions of Bose-Chowla and Singer. The first sentence of Theorem~\ref{Main Theorem} is Observation \#1 in~\cite{GT}.
\begin{theorem}\label{Main Theorem}
    If $h,q,b$ are in the domain of $\Bose$, then $\Bose_h(q,b)$ is a $B_h$ set in $\ZZ/(q^h-1)\ZZ$ with $q$ distinct elements.

    If $h,q,b$ are in the domain of $\Singer$, then $\Singer_h(q,b)$ is a $B_h$ set in $\ZZ/(\tfrac{q^{h+1}-1}{q-1}\ZZ)$ with $q+1$ distinct elements.
\end{theorem}

We say that sets $\cA_1,\cA_2 \in \ZZ/M\ZZ$ are affinely equivalent, writing $\cA_1\sim\cA_2$, if there is some $d$ relatively prime to $M$ and some $s$ with $\cA_2=\{da+s : a\in \cA_1\}$. Clearly, if $\cA_1$ is a $B_h$ set in $\ZZ/M\ZZ$ and $\cA_1\sim\cA_2$, then $\cA_2$ is also a $B_h$ set in $\ZZ/M\ZZ$. We use the notation $d\ast \cA := \{da : a \in \cA\}$ and $\cA+s := \{a+s : a \in \cA\}$ to denote dilations and translations of sets.

We identify when the $B_h$-sets given in Theorem~\ref{Main Theorem} are affinely equivalent in the following theorem.

\begin{theorem}\label{b values}
    Suppose that $h,q=p^r,b$ and $h,q,e$ are in the domain of $\Bose$. If
    \begin{enumerate}[noitemsep,nolistsep,label=(\roman*),itemindent=10mm,labelsep=4mm]
        \item \label{i} $b\equiv e \pmod{\frac{q^h-1}{q-1}}$, or
        \item $\theta^b-\theta^e \in \FF_q$, or
        \item $b\equiv p^i e \pmod{q^h-1}$ for some integer $i$,
    \end{enumerate}
    then $\Bose_h(q,b)\sim \Bose_h(q,e)$.

    Suppose that $h,q,b$ and $h,q,e$ are in the domain of $\Singer$. If
    \begin{enumerate}[resume,noitemsep,nolistsep,label=(\roman*),itemindent=10mm,labelsep=4mm]
        \item $b\equiv e \pmod{\frac{q^{h+1}-1}{q-1}}$, or
        \item $b\equiv p^i e \pmod{q^{h+1}-1}$ for some integer $i$, where $p$ is the prime that divides $q$, or
        \item $\theta^b-\theta^e \in \FF_q$, or
        \item $r\theta^b+t\theta^{e+b}+w\theta^e\in \FF_q$ for some $r,t,w\in \FF_q$, and $r,t$ are not both $0$,
    \end{enumerate}
    then $Singer_h(q,b)\sim \Singer_h(q,e)$.
\end{theorem}

One sad consequence is that for $h=2$ and any allowed $q,b,e$, we have $\Bose_2(q,b)\sim\Bose_2(q,e)$. That is, we do not produce any new (up to affine equivalence) Sidon sets. However,
    \(\Bose_3(5,1) \not\sim \Bose_3(5,4)\)
and for $h>2$ and most (but not all) prime powers $q$ we generate previously unknown $B_h$-sets.

Further, OEIS sequence A227358~\cite{OEIS} gives the minimum diameter of $B_3$-sets with up to 10 elements. We are not aware of any such computation for $h>3$. As comparison, we have also computed the smallest diameters achievable by any subset of any shift of dilations of $\Bose_3(q,b)$, $\Bose_4(q,b)$, $\Singer_3(q,b)$, and $\Singer_4(q,b)$ for all $b$ and small $q$ (projected from $\ZZ/M\ZZ$ to $\ZZ$ in the obvious way). In my opinion, this data suggests that the $\Bose$ and $\Singer$ constructions are not close to optimal for $h>2$, in contrast to the apparent situation for $h=2$.

Let $R_h(n)$ be the maximum possible cardinality of a $B_h$-set contained in $[1,n]$, and let $R_h^{-1}(k)$ be the smallest $n$ such that there is a $B_h$-set with $k$ elements contained in $[1,n]$. The lower bound on $R_h(n)$ and upper bound on $R_h^{-1}(k)$ implied by our constructions is not better than that achieved by Singer's construction alone. Nevertheless, we give several statements using up-to-date results on the distribution of primes, as these results are frequently misstated in the literature.

\begin{theorem}
    For all $n\in\ZZ$, we have $R_h(n) \ge n^{1/h} - 2^{44}n^{154/(155h)}$ and $R_h^{-1}(k) \le k^h + 3^{155h} k^{h-1/155}$.

    If $k,n \ge e^{e^{34}}$, we have $R_h(n) \ge n^{1/h}- 7n^{2/(3h)}$ and $R_h^{-1}(k) \le k^h + (3k)^{h-1/3}$.

    If $k,n$ are sufficiently large, then $R_h(n) \ge n^{1/h}-n^{21/(40h)}$ and $R_h^{-1}(k) \le k^h + 2^h k^{h-19/40}$.

    If the Riemann Hypothesis holds, then
    \[R_h^{-1}(k) < k^h+\log(20k)k^{h-1/2}+2k^{h-1}\log^{2h}(20k),\quad R_h(n) \ge n^{1/h}-(7+\tfrac{\log n}{h})n^{1/(2h)}.\]
\end{theorem}

\section{Two Explicit Examples}
\paragraph{A $\Bose$ Example.} Let $h=3$ and $q=11$; we first compute the various $\Bose_3(11,b)$, and will then give $\Singer_3(11,b)$.

The Conway polynomial for $11^3$ is
    \(C_{11,3}(\theta) = 9+2\theta+\theta^3 \in \FF_{11}[\theta].\)
We have $\FF_{q^3} \cong \FF_q[\theta]/C_{11,3}(\theta)$, and $\theta$ generates the multiplicative group.

Our first task is to find a suitably small set of candidates for $b$. From Theorem~\ref{b values}(i), we only need to consider values between $1$ and $\tfrac{q^h-1}{q-1}=133$, inclusive. As the $\FF_{q^3}$ has only $\FF_q$ as a subfield, only $b=133$ has $\theta^b$ having algebraic degree less than $3$. Further, by Theorem~\ref{b values}(ii) each $b$ is equivalent to $11b$ and $11^2b$. These equivalences combine to give additional equivalences, e.g., $3 \sim 11^2\cdot 3 \sim 97$. The second condition given in Theorem~\ref{b values} is much harder to use, as it requires arithmetic inside the field. For instance,
    \[\theta^{21}-\theta^1 = (\theta^3)^7-\theta=(-9-2\theta)^7-\theta=2^7(1-\theta)^7-\theta=\cdots = 3 \in \FF_{11},\]
and so $1 \sim 21$. With some computerized labor, we find that each $b$ value is equivalent to one of $1$, $2$, $4$, $6$. We have used the Conway polynomial representation, but any finite field representation will lead to four equivalence classes for $b$, but not necessarily these as the smallest representatives of each class.

We then compute inside the field using Definition~\ref{def:Bose} that
    \begin{align*}
        \Bose_3(11,1) &= \{1,21,65,100,111,238,324,523,535,1137,1214\}\\
        \Bose_3(11,2) &= \{2,16,132,237,330,338,389,419,764,1174,1254\}\\
        \Bose_3(11,4) &= \{4,56,116,174,354,626,782,905,979,1147,1183\}\\
        \Bose_3(11,6) &= \{6,152,261,295,311,352,367,891,1092,1113,1228\}
    \end{align*}
By Theorem~\ref{Main Theorem}, these four sets are $B_3$-sets in $\ZZ/1330\ZZ$, and by direct computation we can verify that no two are affinely equivalent. We are not aware of any affine equivalences that are not dictated by Theorem~\ref{b values}.

By directly examining all sets affinely equivalent to these, we notice that
    \[167 \ast \Bose_3(11,6) +330 = \{1,2,27,167,385,397,439,444,484,586,594\}\]
has a particularly small diameter. Consequently $R_3(594) \ge 11$ and $R_3^{-1}(11)\le 594$.

\paragraph{A $\Singer$ Example.} We now compute $\Singer_3(11,b)$. The Conway polynomial for $11^4$ is
\(C_{11,4}(\theta) = 2+10\theta+8\theta^2+\theta^4 \in \FF_{11}[\theta].\)

Our first task is find a suitably small set of candidates for $b$. From Theorem~\ref{b values}(iv), we only need to consider values between $1$ and $\frac{q^{h+1}-1}{q-1}=1464$. We require $\theta^b$ to have algebraic degree $h+1=4$, and that reduces the number of $b$ values to $1452$. Including Theorem~\ref{b values}(v) reduces the number of possible inequivalent $b$ values to $366$. Theorem~\ref{b values}(vi) is significantly more computationally intensive, but reduces the number of inequivalent to $b$ values to at most $36$. Using Theorem~\ref{b values}(vii) is \emph{much} more time-consuming. With the additional assumptions that $r=0,t=1$, we find that each $b$ is equivalent to one of $1$, $2$, $3$, $6$, $8$ or $14$. We have the $B_3$-sets in $\ZZ/1464\ZZ$:
    \begin{align*}
        \Singer_3(11,1) &= \{1,418,502,679,846,1050,1164,1187,1285,1319,1339,1464\}   \\
        \Singer_3(11,2) &= \{2,273,377,432,500,665,674,887,908,1192,1257,1464\} \\
        \Singer_3(11,3) &= \{3, 201, 309, 425, 664, 700, 876, 1061, 1105, 1239, 1357, 1464\} \\
        \Singer_3(11,6) &= \{6,76,388,593,702,734,950,1147,1208,1440,1457,1464\}\\
        \Singer_3(11,8) &= \{8,128,582,624,739,774,841,922,1143,1311,1369,1464\}\\
        \Singer_3(11,14) &= \{14,40,85,492,529,621,683,722,940,969,1151,1464\}
    \end{align*}
By direct computation, no two of these are affinely equivalent. We are not aware of any affine equivalences that are not dictated by Theorem~\ref{b values}.

We further find, after some computation, that
    \[ 481 \ast \Singer_3(11,1) + 102 = \{1,4,36,72,89,102,229,379,583,592,629,738\}.\]
Thus, $R_3(738) \ge 12$ and $R_3^{-1}(12) \le 738$. Moreover,
    \[ 653 \ast \Singer_3(11, 2) + 564 = \{1,22,31,81,92,108,225,406,564,568,592,793\}.\]
Thus, dropping the last element, we find that $R_3(592) \ge 11$ and $R_3^{-1}(11) \le 592$. This is slightly better than the bound from $\Bose_h(11,b)$ sets.

\section{Generalized Bose-Chowla Sets}
Fix an integer $h\geq 2$ and a prime power $q$, and set $M:=q^h-1$. Let $\tau$ be a multiplicative generator of $\FF_{q^h}^\times$ (not necessarily in line with the Conway polynomial). Take $\beta\in \FF_{q^h}$ with algebraic degree $h$.
We define $S_h$ as follows:
    \[S_h(\tau,\beta) := \left\{a \in \ZZ/M\ZZ : \quad \tau^a = \beta +v,\quad v\in \FF_q\right\} .\]
Further, let $\alpha_1,\alpha_2,\dots,\alpha_h$ be a basis for $\FF_{q^h}$ over $\FF_q$ as a vector space, with $\alpha_1=1$ and $\alpha_2=\beta$.

As $1,\beta,\dots,\alpha_h$ is a basis, each $x' \in  \FF_q$ corresponds to a distinct $x \in \ZZ/M\ZZ$ by the equation $\theta^x=1\cdot x'+1\cdot\beta+\sum_{i=3}^h 0\cdot \alpha_i$, so that $S_h(\tau,\beta)$ has exactly $q$ elements.

Consider $k \in \ZZ/M\ZZ$, and suppose  that $a_1,\dots,a_h,b_1,\dots,b_h \in S_h(\tau,\beta)$ satisfy
$$k = a_1+\dots+a_h \equiv b_1+\dots+b_h  \pmod{M}.$$
As $\tau$ has multiplicative order $q^h-1=M$, we have
\begin{equation*}
    \tau^k
    = \tau^{a_1+\dots+a_h}
    =\prod_{i=1}^h \tau^{a_i}
    = \prod_{i=1}^h (\beta + a_i')
\end{equation*}
for some $a_i'\in\FF_q$
In the same manner,
    \[\tau^k = \prod_{i=1}^h (\beta + b_i').\]

Now define polynomials $f,g\in\FF_q[x]$ by
    \[f(x) = \prod_{i=1}^h (x+a_i'),\qquad g(x)  = \prod_{i=1}^h (x+b_i').\]
Then $\beta$ (which has algebraic degree $h$) is a root of $f(x)-g(x)$ (which has degree at most $h-1$), from which we learn that $f(x)-g(x)$ is identically $0$, i.e., $f(x)=g(x)$. We have unique factorization over finite fields, so that
    \[\left\{ a_1',\dots,a_h' \right\} = \left\{b_1',\dots,b_h'\right\}\]
as multisets. As noted above, that $\alpha_1,\dots,\alpha_h$ is a basis implies that $a_i',b_i'\in\FF_q$ uniquely define $a_i,b_i$ in $\ZZ/M\ZZ$, and consequently
    \[\left\{ a_1,\dots,a_h \right\} = \left\{b_1,\dots,b_h\right\}\]
as multisets. That is, $S_h(\tau,\beta)$ is a $B_h$-set  in $\ZZ/M\ZZ$.

We can take $\tau$ to be $\theta$, the generator provided in the Conway polynomial representation of $\FF_{q^h}$, and note that $\beta=\theta^b$ for some $b\in\ZZ/M\ZZ$, so that $S_h(\tau,\beta)=\Bose_h(\theta,b)$. We have proven the claims in Theorem~\ref{Main Theorem} concerning $\Bose_h(q,b)$ sets.

Before proceeding into the proof of Theorem~\ref{b values} as it pertains to $\Bose$ sets, we spend a few words noting some tempting generalizations that aren't really meaningful generalizations. First, fix any basis $\alpha_1,\dots,\alpha_h$ of $\FF_{q^h}$ over $\FF_q$, and any constants $c_1,\dots,c_{h-1}\in\FF_q$, not all $0$ and with $(c_1\alpha_1 + \dots + c_{h-1}\alpha_{h-1})\alpha_h^{-1}$ having degree $h$. Then, the set
    \[\{a \in \ZZ/M\ZZ \colon  \tau^a=c_1\alpha_1+\dots+c_{h-1}\alpha_{h-1}+v\alpha_h,v\in\FF_q\}\]
is a $B_h$-set with $q$ elements. By details we spare the reader, each such set is affinely equivalent to $\Bose_h(q,b)$ for some integer $b$. Second, we note that if $\tau$ is also a generator of the multiplicative group of $\FF_{q^h}$, then $S_h(\tau,\beta) \sim S_h(\theta,\beta).$ Specifically, $\tau=\theta^t$ for some $t$, and since $\tau$ is a generator, $\gcd(t,M)=1$; let $t^{-1}$ be the inverse of $t$ modulo $m$. Then
    \begin{multline*}
        S_h(\tau,\beta)
            := \{ a \in \ZZ/M\ZZ : \tau^a = \beta + v, \quad v \in \FF_q\} \\
            = \{ a \in \ZZ/M\ZZ : \theta^{ta} = \beta + v, \quad v \in \FF_q\}
            = t^{-1} \ast S_h(\theta,\beta).
    \end{multline*}

We now turn to the task of determining when
\[\Bose_h(q,b) \sim \Bose_h(q,e).\]

First, suppose that $b\equiv e \pmod{\tfrac{q^h-1}{q-1}}$. Then for some integer $x$ we have $b=e+x\tfrac{q^h-1}{q-1}$ and $\theta^b=\theta^e\theta^{x(q^h-1)/(q-1)}=w\theta^e$, and $w=(\theta^{(q^h-1)/(q-1)})^x \in\FF_q$ because $\theta^{(q^h-1)/(q-1)}$ is in $\FF_q$. We have
\begin{align*}
    \Bose_h(q,b)
    &:= \left\{a \in \ZZ/M\ZZ : \quad \theta^a = \theta^b +v,\quad v\in \FF_q\right\} \\
    &= \left\{a \in \ZZ/M\ZZ : \quad \theta^a = w\theta^e +v,\quad v\in \FF_q\right\} \\
    &= \left\{a \in \ZZ/M\ZZ : \quad \theta^{a-x(q^h-1)/(q-1)} = \theta^e +vw^{-1},\quad v\in \FF_q\right\} \\
    &= x\tfrac{q^h-1}{q-1}+\left\{a \in \ZZ/M\ZZ : \quad \theta^a = \theta^e +v,\quad v\in \FF_q\right\} \\
    &=  x\tfrac{q^h-1}{q-1}+\Bose_h(q,e).
\end{align*}
Thus, $\Bose_h(q,b) \sim \Bose_h(q,e)$.

Now, suppose that $pb \equiv e \pmod{M}$, where $p$ is the characteristic of the field $\FF_{q^h}$. The map $u \mapsto u^p$, the Frobenius automorphism, is a bijection and satisfies $(u+v)^p =u^p+v^p$ for any $u,v\in \FF_{q^h}$. We have
\begin{align*}
    \Bose_h(q,b)
    &:= \left\{a \in \ZZ/M\ZZ : \quad \theta^a = \theta^b +v,\quad v\in \FF_q\right\} \\
    &= \left\{a \in \ZZ/M\ZZ : \quad (\theta^a)^p = (\theta^b +v)^p,\quad v\in \FF_q\right\} \\
    &= \left\{a \in \ZZ/M\ZZ : \quad \theta^{ap} = \theta^{pb} +v^p,\quad v\in \FF_q\right\} \\
    &= p^{-1}\ast\left\{a \in \ZZ/M\ZZ : \quad \theta^{ap} = \theta^{e} +v,\quad v\in \FF_q\right\} \\
    &=p^{-1} \ast \Bose_h(q,e).
\end{align*}
It follows that if $b \equiv p^i e \pmod{M}$ for any $i$, then $\Bose_h(q,b)\equiv\Bose_h(q,e)$.

Now suppose  that  $w:=\theta^e-\theta^b \in \FF_q$. Then
\begin{align*}
    \Bose_h(q,b)
    &:= \left\{a \in \ZZ/M\ZZ : \quad \theta^a = \theta^b +v,\quad v\in \FF_q\right\} \\
    &= \left\{a \in \ZZ/M\ZZ : \quad \theta^a = \theta^e-w +v,\quad v\in \FF_q\right\} \\
    &= \Bose_h(q,e).
\end{align*}
Thus, $S_h(q,\theta,\theta^b) \sim S_h(q,\theta,\theta^e)$.

This concludes the proof of all of the claims regarding $\Bose$ sets made in Theorems~\ref{Main Theorem} and~\ref{b values}.

\section{Generalized Singer Sets}

Fix an integer $h\geq 2$ and a prime power $q$, and set $M:=(q^{h+1}-1)/(q-1)$. Let $\tau$ be a multiplicative generator $\FF_{q^{h+1}}^\times$. Suppose further that $\beta$ has algebraic degree $h+1$. We define $S_h$ as follows:
    \[S_h(\tau,\beta) := \left\{ a\in \ZZ/M\ZZ :\quad
        \tau^a = u \beta +v, \quad u,v \in \FF_q\right\}.\]
Further, let $\alpha_1,\alpha_2,\dots,\alpha_h$ be a basis for $\FF_{q^h}$ over $\FF_q$ as a vector space, with $\alpha_1=1$ and $\alpha_2=\beta$. Note that $\tau^M \in \FF_q$, as is $\tau^{kM}$ for any integer $k$.

We first argue that $S_h(\tau,\beta)$ has $q+1$ distinct elements. As $1,\beta,\alpha_3,\dots,\alpha_h$ is a basis, for each $u,v\in\FF_q$, not both $0$, there is a unique $a$ in $[1,q^{h+1}-1]$ with $\tau^a=u+v\beta$. That is, there are $q^2-1$ such $a$. For each particular $a$, there is also a solution (with different $u,v$) with $a+kM$ for any integer $k$, as
    \[ \tau^{a+kM}=\tau^{kM}\big( u\beta +v \big)= (wu)\beta+(wv),\]
where $w=\tau^{kM}\in \FF_q$, and  so  $wu,wv\in\FF_q$. Thus, the $q^2-1$ solutions with $1\leq  a \leq q^{h+1}-1$ fall into congruence classes modulo $M$. Each congruence class modulo $M$ has $q-1$ elements in $1\leq a \leq q^{h+1}-1$, so that $S_h(\tau,\beta)$ consists of $(q^2-1)/(q-1)=q+1$ distinct elements.

We now prove that $S_h(\tau,\beta)$ is a $B_h$-set in $\ZZ/M\ZZ$. Define functions $u,v\colon S_h(\tau,\beta) \to \FF_q$ by
    \[\tau^k = u(k) \beta + v(k).\]
Note that $u(M)=0$ since $\tau^M\in\FF_q$.
Clearly the pair $(u(k),v(k))$ uniquely determines $k\in S_h(\tau,\beta)$. But more surprisingly, the value $-v(k)u(k)^{-1}$ (possibly undefined) uniquely determines $k\in S_h(\tau,\beta)$, as we now explain. Suppose $-v(k)u(k)^{-1}$ is undefined, whence $u(k)=0$. Then $\tau^k=v(k)\in \FF_q$, so that $k\equiv 0 \pmod{M}$. Since $S_h(\tau,\beta) \subseteq \ZZ/M/ZZ$, we must have $k=M$. Otherwise, $w:=-v(k)u(k)^{-1}$ is defined, whence $w u(k)=v(k)$. This means that $\tau^k = u(k) \beta + wu(k)$. We know $\beta$ and $w$, and therefore know the value $y\in[1,q^{h+1}-1]$ with $\tau^y = \beta+w$. Since $u(k)\in\FF_q$, we have $y\equiv k \pmod{M}$, whence there is a unique candidate for $k$ in $[1,M]$.

Now suppose that
\begin{equation}\label{eq:2h terms}
    a_1+\dots+a_h \equiv b_1+\dots+b_h \pmod M,
\end{equation}
with $a_i,b_i \in S_h(\tau,\beta)$, and we must show that
\begin{equation}\label{eq:goal}
    \{a_1,\dots,a_h\} = \{b_1,\dots,b_h\}
\end{equation}
as multisets.
From line \eqref{eq:2h terms}, there is an integer $x$ with
    \[a_1+\dots+a_h = kM+ b_1+\dots+b_h.\]
Let $w=\tau^{kM}\in\FF_q$. We then have, using that $a_i,b_i\in S_h(\tau,\beta)$,
\begin{multline*}
    \prod_{i=1}^h \big(u(a_i) \beta + v(a_i) \big)
    = \prod_{i=1}^h \tau^{a_i}
    = \tau^{\sum_{i=1}^h a_i}
    = \tau^{xM+\sum_{i=1}^h b_i} \\
    = w \prod_{i=1}^h \tau^{b_i}
    = w \prod_{i=1}^h \big(u(b_i) \beta + v(b_i) \big)
\end{multline*}
We define the polynomials (each  with degree at most $h$) in $\FF_q[x]$
\begin{equation*}
    f(x) := \prod_{i=1}^h \big(u(a_i) x + v(a_i) \big), \qquad
    g(x) := \prod_{i=1}^h \big(u(b_i) x + v(b_i) \big).
\end{equation*}
Then $\beta$, which by  hypothesis  has degree $h+1$, is  a root  of the polynomial $f(x)-wg(x)$, which has  degree at most $h$. Thus  $f(x)=w  g(x)$, and $f,g$ must have the same roots in the same multiplicities. That is, the multisets are equal
    \[\left\{ -v(a_i)u(a_i)^{-1} \colon 1\leq i \leq h \right\}
        =\left\{ -v(b_i)u(b_i)^{-1} \colon 1\leq i \leq h \right\},\]
including the number of  occurrences of undefined elements. As noted above, the value of $-v(a_i)u(a_i)^{-1}$ uniquely determines $a_i$, so that the multiset equality
    \[\{ a_1, \dots,a_h\} = \{b_1,\dots,b_h\}\]
holds.

We can take $\tau$ to be $\theta$, the generator provided in the Conway polynomial representation, and we can locate $b\in\ZZ/M\ZZ$ so that $\beta=\theta^b$, and then $S_h(\tau,\beta)=\Singer_h(\theta,b)$. We have proven the claims in Theorem~\ref{Main Theorem} concerning $\Singer_h(q,b)$ sets.

Before proceeding into the proof of Theorem~\ref{b values} as it pertains to $\Singer$ sets, we note that as with $\Bose$ sets, neither the completion of $1,\beta$ into a basis (which we do not elaborate upon) nor the particular choice of generator actually matters, up to affine equivalence, which we now elaborate. Suppose $\tau=\theta^k, \beta=\theta^b$. Then
\begin{align*}
    S_h(\tau,\beta)
    &:= \left\{a \in\ZZ/M\ZZ :
        \quad\tau^a = u \beta +v, \quad u,v \in \FF_q\right\} \\
    &= \left\{a \in\ZZ/M\ZZ :
        \quad\theta^{ka} = u \theta^{b} +v, \quad u,v \in \FF_q\right\}
    &= k^{-1}\ast \Singer_h(q,b).
\end{align*}
Thus, we have lost nothing by defining $\Singer_h(q,b)$ with respect to a specific generator.

Now suppose that $b\equiv e \pmod{M}$,  whence $b=e+kM$ for some integer $k$ and
$$\theta^{b} = \theta^e\theta^{kM}= w \theta^e$$
for some $0\neq w\in \FF_q$. We have
\begin{align*}
    \Singer_h(q,b)
    &:= \left\{a \in\ZZ/M\ZZ :
        \quad\theta^a = u \theta^b +v, \quad u,v \in \FF_q\right\} \\
    &= \left\{a \in\ZZ/M\ZZ :
        \quad\theta^a = u w \theta^e +v, \quad u,v \in \FF_q\right\}
    &= \Singer_h(q,e).
\end{align*}

Recall that the Frobenius automorphism $u\mapsto u^p$, where $p$ is the characteristic of $\FF_{q^h}$ fixes each element of $\FF_q$, and satisfies the ``children's binomial theorem'': $(u+v)^p = u^p + v^p$ for all $u,v\in\FF_{q^h}$. Suppose that $e\equiv p b \pmod{M}$. Then
\begin{align*}
    \Singer_h(q,b)
    &:= \left\{a \in\ZZ/M\ZZ :
        \quad\theta^a = u \theta^b +v, \quad u,v \in \FF_q\right\} \\
    &= \left\{a \in\ZZ/M\ZZ :
        \quad\theta^{ap} = u^p \theta^{bp} + v^p, \quad  u,v \in \FF_q\right\} \\
    &= \left\{a \in\ZZ/M\ZZ :
    \quad\theta^{ap} = u \theta^{bp} + v, \quad  u,v \in \FF_q\right\} \\
    &= p^{-1} \ast \Singer_h(q,bp).
\end{align*}
It follows that if $b\equiv p^i e \pmod{q^{h+1}-1}$ for some integer $i$, then $\Singer_h(q,b)\sim \Singer_h(q,e)$.

While Theorem~\ref{b values}(vi) is a special case of Theorem~\ref{b values}(vii), we provide a separate proof of the easier (vi) as it is independently useful in computations.
Suppose that $w:=\theta^e-\theta^b \in \FF_q$. Then
\begin{align*}
    \Singer_h(q,b)
    &:= \left\{a \in\ZZ/M\ZZ :
    \quad\theta^a = u \theta^b +v, \quad u,v \in \FF_q\right\} \\
    &=  \left\{a \in\ZZ/M\ZZ :
    \quad\theta^a = u (\theta^e-w) +v, \quad u,v \in \FF_q\right\} \\
    &=  \left\{a \in\ZZ/M\ZZ :
    \quad\theta^a = u \theta^e +v-uw, \quad u,v \in \FF_q\right\} \\
    &=  \Singer_h(q,e).
\end{align*}

We now address Theorem~\ref{b values}(vii). Suppose that $r,t,w,s \in \FF_q$, and at least one of $r,t$ is nonzero, and
    \[r \theta^b + t\theta^{e+b}+w\theta^e = s.\]
Then $(r+t\theta^e)\theta^b = s-w\theta^e$. Since $1,\theta^e$ are linearly independent over $\FF_q$ and at least one of $r,t$ is nonzero, we know that $r+t\theta^e$ is nonzero, say $\theta^k=r+t\theta^e$. We have
\begin{align*}
    \Singer_h(q,b)
    &:= \left\{a \in\ZZ/M\ZZ :
        \quad\theta^a = u \theta^b +v, \quad u,v \in \FF_q\right\} \\
    &=  \left\{a \in\ZZ/M\ZZ :
        \quad\theta^k \theta^a = (r+t\theta^e)(u \theta^b +v), \quad u,v \in \FF_q\right\} \\
    &= \left\{a \in\ZZ/M\ZZ :
        \quad\theta^k \theta^a = rv +u(r+t\theta^e)\theta^b+rt\theta^e, \quad u,v \in \FF_q\right\} \\
    &= \left\{a \in\ZZ/M\ZZ :
        \quad\theta^k \theta^a = rv +u(s-w\theta^e)+rt\theta^e, \quad u,v \in \FF_q\right\} \\
    &= \left\{a \in\ZZ/M\ZZ :
        \quad\theta^k \theta^a = rv +us+(rt-uw)\theta^e, \quad u,v \in \FF_q\right\} \\
    &= -k + \left\{a \in\ZZ/M\ZZ :\quad \theta^a = u' \theta^e + v', \quad u',v'\in \FF_q\right\}. \\
        &= -k + \Singer_h(q,e)
\end{align*}
The equality of the last line relies upon the nonsingularness of the equations $u'=rv+us,v'=rt-uw$, which follows from $\theta^b$ being outside $\FF_q$ and the equation $(r+t\theta^e)\theta_b = s-w\theta^e$

This concludes the proof of all of the claims regarding $\Singer$ sets made in Theorems~\ref{Main Theorem} and~\ref{b values}.

\section{Lower Bounds on $\mathbf{R_h(n)}$ and upper bounds on $\mathbf{R_h^{-1}(k)}$}
Computing $R_2^{-1}(k)$ is an old game already~\cite{1906.Veblen}. Babcock~\cite{1953.Babcock} computed by hand for $k\le 10$ (his value of $R_2^{-1}(10)$ is incorrect). More recently, the OGR Project~\cite{OGR} has computed $R_2^{-1}(28)=585$; the computation took $8.5$ years on thousands of machines. We refer the reader to the Wikipedia page for Golomb Rulers~\cite{wikiGolomb} for $R_2^{-1}(k)$ for $k\le 28$ and for the sets that are optimal.

Another massive computation for $R_2^{-1}(k)$ was carried out by Dogon \& Rokicki~\cite{2016.Dogon.Rokicki}. With several clever optimizations, they computed the bound achieved by all subsets of all sets affinely equivalent to $\Bose_2(q,1)$ and $\Singer_2(q,1)$ for all $q\le 40\,000$. In this section, we report on a similar computation, much smaller in scale, for $h=3$ and $h=4$.

The asymptotic growth of $R_h(n)$ (respectively, $R_h^{-1}(k)$) is not known for $h>2$. The best lower bounds on $R_h(n)$ (upper bounds on $R_h^{-1}(k)$) arise from the construction of Singer~\cite{Singer}. Our generalization produces many more such sets, but they are of roughly the same size. Nevertheless, we feel it would be a contribution to the literature to record the resulting bounds under several hypotheses.

Clearly, $R_h^{-1}(1)=1$, $R_h^{-1}(2)=2$, and $R_h^{-1}(3)=\max\{1,2,h+2\}=h+2$. We therefore restrict our attention to $k\ge 4$ and $n\ge h+3$.

By ``Bose $B_h$-set'', we mean any affine image of $\Bose_h(q,\theta,b)$ for any $q,\theta,b$ in the domain of $\Bose_h$.     By ``Singer $B_h$-set'', we mean any affine image of $\Singer_h(q,b)$ for any $q,\theta,b$ in the domain of $\Singer_h$.

While Singer $B_h$-sets are slightly thicker than Bose $B_h$-sets, it is easier to work with Bose sets. First, if $q$ is a prime power, then $\Bose_h(q,1)$ is a set with $q$ elements modulo $q^h-1$. Thus, $R_h(q^h-1) \geq q$ and $R_h^{-1}(q) \le q^h-1$. Thus, if $k \le q$ then $R_h^{-1}(k) \le R_h^{-1}(q) < q^h$. The is now reduced to locating a prime power greater than $k$, but not too much greater.

We will state results that work for every $k$, for $k>k_0$ with explicit $k_0$, and for $k$ sufficiently large assuming the Riemann Hypothesis. The bounds are either impracticably bad for small $k$, or only apply for impracticably large $k$, or use an impracticably difficult hypothesis. Except for $h=2$, we do not believe that the main terms reported below even have the ``correct'' coefficient.
We start with the most explicit unconditional result.
\begin{theorem*}[Cully-Hugill~\cite{Cully-Hugill.155}]
    For all integers $n\ge 1$, there is a prime between $n^{155}$ and $(n+1)^{155}$.
\end{theorem*}
It follows that there is a prime between $\lceil k^{1/155} \rceil^{155}$ and $\lceil k^{1/155}+1\rceil^{155}$. As
    \[\lceil k^{1/155}+1\rceil^{155h} < (k^{1/155}+2)^{155h} < k^h+3^{155h}k^{h-1/155},\]
we have the statement in the theorem below for $R_h^{-1}(k)$. Assuming that
    \[k^{155}<q<(k+1)^{155}<n^{1/h} \le (k+2)^{155}\]
and using the straightforward $k^{155}>(k+2)^{155}-2^{44}(k+1)^{154}$ yields the $R_h(n)$ result.
\begin{theorem}
    For all $k\ge 4$ and $n\ge h+3\ge 5$, we have $R_h^{-1}(k) < k^h+3^{155h} k^{h-1/155}$ and $R_h(n) \ge n^{1/h}-2^{44}n^{154/(155h)}$.
\end{theorem}

%
%

For large $k$, we can do somewhat better.
\begin{theorem*}[Cully-Hugill~\cite{Cully-Hugill.155}]
    For all integers $n>\exp(\exp(32.537))$, there is a prime between $n^3$ and $(n+1)^3$.
\end{theorem*}
Hence:
\begin{theorem}
    For all $k>e^{e^{34}}$, we have $R_h^{-1}(k) < k^{h}+(3k)^{h-1/3}$ and $R_h(n) > n^{1/h}-7n^{2/(3h)}$.
\end{theorem}

The following famed result~\cite{MR1851081} is beautifully straightforward to use.
\begin{theorem*}[Baker \& Harman \& Pintz~\cite{MR1851081}]
    If $x$ is sufficiently large, then there is a prime in the interval $[x-x^{21/40},n]$, and in the interval $[x,x+x^{21/40}]$.
\end{theorem*}
This leads to:
\begin{theorem}
    If $k,n$ are sufficiently large, then $R_h^{-1}(k) < k^h+2^h k^{h-19/40}$ and $R_h(n) \ge n^{1/h}-n^{21/(40h)}$.
\end{theorem}

Assuming the Riemann Hypothesis, we naturally have stronger results. The best result along these lines of which this author is aware follows~\cite{MR3498632}.
\begin{theorem*}[Dudek \& Greni\'{e} \& Lo\"{\i}c~\cite{MR3498632}]
    Assuming the Riemann Hypothesis, for all $n\ge 2$, there is a prime between $n^2$ and $(n+(1+\tfrac{1}{\log n})^2\log n)^2$.
\end{theorem*}
This leads directly to the following.
\begin{theorem}
    Assume the Riemann Hypothesis, and that $k\ge 4$, $n\ge h+3$. Then
        \[R_h^{-1}(k) < k^h+\log(20k)k^{h-1/2}+2k^{h-1}\log^{2h}(20k),\quad R_h(n) \ge n^{1/h}-(7+\tfrac{\log n}{h})n^{1/(2h)}.\]
\end{theorem}

\section{Explicit Computations}
For $k\le 9$, we computed the minimum-diameter $B_3$-sets in $\ZZ$ by brute force. This allowed us to find the sequence in the OEIS (A227358), where $R_3^{-1}(10)$ is also reported. These results are shown in Table~\ref{tab:optimal B3}.
\begin{table}[p]
    \[\begin{array}{ccc}
        k & R_3^{-1}(k) & \text{witness} \\ \hline
        1 & 1 & \{0\} \\
        2 & 2 & \{0,1\} \\
        3 & 5 & \{0,1,4\} \\
        4 & 12 & \{0, 1, 7, 11\}, \{0, 1, 8, 11\}\\
        5 & 24 & \{0, 1, 15, 18, 23\}, \{0, 1, 15, 20, 23\}\\
        6 & 46 & \{0, 2, 11, 26, 42, 45\}\\
        7 & 83 & \{0, 1, 7, 50, 59, 78, 82\}, \{0, 2, 23, 45, 72, 79, 82\}\\
        & & \{0, 4, 23, 32, 75, 76, 82\}\\
        8 & 130 & \{0, 2, 5, 34, 74, 107, 120, 129\}\\
        9 & 209 & \{0, 1, 17, 26, 127, 138, 185, 204, 208\}\\
        & & \{0, 1, 18, 76, 83, 162, 188, 193, 208\} \\
        10 & 310 & \\
    \end{array}\]
    \caption{A227358, computations by John Tromo, sets and $k\le 9$ independently computed by the author.}
    \label{tab:optimal B3}
\end{table}

We have computed all translations of all dilations of all subsets of the Singer and Bose $B_3$-sets generated with small $q$ and any $b$. These results are shown in Table~\ref{tab:B3 from BoseSinger}. The same computation was performed for $B_4$-sets, and those results are given in Table~\ref{tab:B4 from BoseSinger}.
\begin{table}[p]
\[
\begin{array}{r|cccccc}
    k & R_3^{-1}(k) & \text{from Greedy} & \text{from $\Bose$} & \text{with }q & \text{from $\Singer$} & \text{with }q \\ \hline
    1 & 1 & 1& 1 & 2 & 1 & 2 \\
    2 & 2 & 2& 2 & 2 & 2 & 2 \\
    3 & 5 & 5& 5 & 4 & 5 & 2 \\
    4 & 12 & 14& 12 & 5 & 14 & 3 \\
    5 & 24 & 33& 33 & 5 & 28 & 4 \\
    6 & 46 & 72& 73 & 11 & 57 & 5 \\
    7 & 83 & 125& 122 & 7 & 121 & 7 \\
    8 & 130 & 219 & 202 & 8 & 157 & 7 \\
    9 & 209 & 376 & 306 & 9 & 258 & 8 \\
    10 & 310 & 573 & 493 & 11 & 365 & 9 \\
    11 & {} & 745 & 594 & 11 & 592 & 11 \\
    12 & {} & 1209 & 894 & 13 & 738 & 11 \\
    13 & {} & 1557 & 1044 & 13 & 1014 & 13 \\
    14 & {} & 2442 & 1612 & 17 & 1236 & 13 \\
    15 & {} & 3098 & 1874 & 17 & 1877 & 16 \\
    16 & {} & 4048 & 2247 & 16 & 2071 & 16 \\
    17 & {} & 5298 & 2537 & 17 & 2392 & 16 \\
    18 & {} & 6704 & 3433 & 19 & 2960 & 17 \\
    19 & {} & 7839 & 3821 & 19 & 3679 & 19 \\
    20 & {} & 10987 & 5578 & 23 & 4326 & 19 \\
    21 & {} & 12332 & 6060 & 23 & 5849 & 23 \\
    22 & {} & 15465 & 6212 & 23 & 6476 & 23 \\
    23 & {} & 19144 & 6997 & 23 & 7229 & 23 \\
    24 & {} & 24546 & 8846 & 25 & 8010 & 23 \\
    25 & {} & 28974 & 9624 & 25 & 8854 & 25 \\
    26 & {} & 34406 & 11447 & 27 & 10177 & 25 \\
    27 & {} & 37769 & 12088 & 27 & 12143 & 27 \\
    28 & {} & 45864 & 14272 & 29 & 13432 & 27 \\
    29 & {} & 50877 & 15544 & 29 &  &  \\
    30 & {} & 61372 & 17999 & 31 &  &  \\
\end{array}
\]
\caption{The upper bounds on $R_3^{-1}$ that arise from Singer and Bose $B_3$-sets, and also the greedy $B_3$-set (A096772).}
\label{tab:B3 from BoseSinger}
\end{table}

\begin{table}[p]
\[
\begin{array}{r|cccccc}
    k  & R_4^{-1}(k) & \text{from Greedy}  & \text{from $\Bose$} & \text{with }q & \text{from $\Singer$} & \text{with }q \\ \hline
    1  & 1           & 1                   & 1                   & 2             & 1                     & 2             \\
    2  & 2           &  2                  & 2                   & 2             & 2                     & 2             \\
    3  & 6           &   6                 & 6                   & 3             & 6                     & 2             \\
    4  & 16          &    22               & 26                  & 5             & 18                    & 3             \\
    5  & 42          &      56             & 89                  & 5             & 71                    & 5             \\
    6  & 101         &        154          & 212                 & 7             & 156                   & 5             \\
    7  &             &           369       & 404                 & 7             & 388                   & 7             \\
    8  &             &              857    & 959                 & 8             & 693                   & 7             \\
    9  &             &                1425 & 1731                & 11            & 1290                  & 9             \\
    10 &             & 2604                & 2878                & 11            & 2345                  & 9             \\
    11 &             &     4968            & 4469                & 11            & 4053                  & 11            \\
    12 & \text{}     &         8195        & 7967                & 13            & 5174                  & 11            \\
    13 & \text{}     &             13664   & 9903                & 13            & 9328                  & 13            \\
    14 & \text{}     &             22433   & 15907               & 16            & 11348                 & 13            \\
    15 & \text{}     & 28170               & 20849               & 16            &                       &               \\
    16 & \text{}     &      47689          & 25397               & 16            &                       &               \\
    17 & \text{}     &           65546     & 35282               & 17            &                       &               \\
    18 & \text{}     &           96616     & 45783               & 19            &                       &               \\
    19 & \text{}     &   146249            & 58033               & 19            &                       &
\end{array}
\]
\caption{The upper bounds on $R_4^{-1}$ that arise from Singer and Bose $B_4$-sets, and the greedy $B_4$-set (A365300).}
\label{tab:B4 from BoseSinger}
\end{table}

\section{Open Questions}
The following questions are interesting to the author, who does not know of solutions.
\begin{enumerate}
    \item The greedy $B_2$-set is called the Mian-Chowla sequence~\cite{1944.Chowla.Mian}, and the first terms were computed in the 1940s. I'm not aware of any computation of the greedy $B_h$ sequence for $h>2$. I have added these sequences to the OEIS for $4\le h\le 9$ (sequences A365300 through A365305).
    \item The conditions in Theorem~\ref{b values} are necessary for $\Bose_h(q,e)\sim\Bose_h(q,b)$; are they sufficient? Also, for Singer sets.
    \item Is there a faster way to interpret Theorem~\ref{b values}(ii)? Theorem~\ref{b values}(vii) is particularly time consuming, can one assume without loss of generality that $r=0$ and $t=1$?
    \item Does $\Bose_2(q,1)$ always have two elements whose difference is relatively prime to $q^2-1$? Equivalently, is there an affine image of $\Bose_2(q,\theta,1)$ that contains $\{1,2\}$? Is there any $m,s,q$ with
        \[\{0,1,4,10,18,23,25\}\subseteq m \ast \Bose_2(q,1)+s \pmod{q^2-1} ? \]
    Halberstam \& Laxton~\cite{1963.Halberstam.Laxton} considered the $m$ for which there is an $s$ with $\Bose_2(q,1) = m\ast \Bose_2(q,1)+s$. Can this be generalized to $h>2$?
    Also for Singer sets.
    \item Does the largest modular gap between consecutive elements  of $\Bose_2(q,1),\Singer_2(q,1)$ have order $O(q)$? It  seems not, even if one chooses an affine image to make the largest gap as small as possible.
    \item It is obvious that affine maps preserve the $B_h$ property. The existence of Bose sets that are not affine images of each other suggests that there may be some more general arithmetic (or geometric) operation (beyond affine equivalence) that preserves the $B_h$ property in cyclic groups.
\end{enumerate}


\appendix

\end{document}